\documentclass[11pt,fleqn]{article}

\usepackage[cp1251]{inputenc}

\usepackage{amsmath,amssymb,amsthm,enumerate,epsfig,graphics,cite}
\usepackage[ukrainian,english]{babel}

\newcommand{\todo}[1][\null]{\ensuremath{\clubsuit}}

\newcommand{\noprint}[1]{}

\begin{document}
\begin{center}
\Large\bf
Logarithmic Integration Method for Solving of First and Second Order Differential Equations
\end{center}
\begin{center}
\it
A. Ponomarenko
\end{center}

{\it

In this article we present logarithmic methods for solving first order and second order ordinary differential equations. The essence of the method is that we apply the basic properties derivatives and logarithms to reduce the number of terms in the equation. Here we carry out this only for equations of the first and second order. Similar methods can also be used to obtain solutions to higher order equations.
}
\medskip


\medskip

The main methods for solving ordinary differential equations have long been studied and are known to everyone. In \cite {1}, \cite {2}, \cite {3} presents some basic methods for integrating simple ordinary differential equations (ODEs) and focuses on real solutions of ODEs with real coefficients. It describes homogeneous linear equations with constant coefficients. In \cite {1} shows that the general solution of nonhomogeneous linear equations with constant coefficients is the sum of the complementary function (the general solution of the corresponding homogeneous equation) and a particular integral. This article discusses a new approach to solving ordinary differential equations using the simplest elementary operations. The calculations can be cumbersome, but we do not lose particular solutions to differential equations.

\medskip

Let $f(x)$, $g(x)$ be Riemann integrable functions; $y=y(x)$, $y'(x)=\frac{dy(x)}{dx}$, $y''(x)=\frac{d^{2}y(x)}{dx^{2}}$, $ \log y=\ln y =\log_{e}y $; $C$,
$C_{1}$, $C_{2}$, $C_{1,1}$, ... , $C_{1,7}$, $C_{2,1}$, $C_{2,2}$ is an integration constant. The symbol $\Rightarrow$ between two formulas will mean that the second formula follows from the first one.

\medskip
\section{First order differential equations}

{\it 1.1. Linear inhomogeneous first order differential equation:}
\begin{gather}\label{eq1}
y'(x)+f(x)y(x)=g(x).
\end{gather}
{\it Logarithmic integration method.} In equation \eqref{eq1} the function $g(x)$ is not identically zero. Then $y(x)$ be not identifically zero. Then with equations \eqref{eq1} we get
$$
\frac{y'(x)}{y(x)}+f(x)=\frac{g(x)}{y(x)}, \quad \Rightarrow
$$
$$
(\log |y(x)|)'+f(x)=\frac{g(x)}{y(x)}, \quad \Rightarrow
$$
\begin{gather}\label{eq1a}
(\log |y(x)|)'+\left(\int f(x)dx\right)'=\frac{g(x)}{y(x)}, \quad \Rightarrow
\end{gather}
$$
(\log |y(x)|)'+\left(\log e^{\int f(x)dx}\right)'=\frac{g(x)}{y(x)}, \quad \Rightarrow
$$
$$
\left(\log |y(x)|+\log e^{\int f(x)dx}\right)'=\frac{g(x)}{y(x)}, \quad \Rightarrow
$$
$$
\left(\log \left(|y(x)| e^{\int f(x)dx}\right)\right)'=\frac{g(x)}{y(x)}, \quad \Rightarrow
$$
$$
\frac{\left(y(x) e^{\int f(x)dx}\right)'}{y(x) e^{\int f(x)dx}}=\frac{g(x)}{y(x)}, \quad \Rightarrow
$$
$$
\left(y(x) e^{\int f(x)dx}\right)' =g(x) e^{\int f(x)dx}, \quad \Rightarrow
$$
$$
y(x) e^{\int f(x)dx} =\int g(x) e^{\int f(x)dx} dx + C, \quad \Rightarrow
$$
\begin{gather}\label{eq2}
y(x)  =e^{-\int f(x)dx} \left( \int g(x) e^{\int f(x)dx} dx + C \right).
\end{gather}

{\it Remark 1.1.1.} A similar method can be obtain solution the equation \eqref{eq1} in the Cauchy form:
\begin{gather}\label{eq2d}
y(x)  =e^{-\int_{x_0}^x f(t)dt}\left(\int_{x_0}^x g(\tau) e^{\int_{x_0}^\tau f(\sigma)d\sigma} d\tau + y(x_{0}) \right),
\end{gather}
where $x_{0}$ is a given constant.
Indeed, the equation \eqref{eq1a} is equivalent to the equation
\begin{gather}\label{eq1b}
(\log |y(x)|)'+\left(\int f(x)dx+ C_{1}\right)'=\frac{g(x)}{y(x)},
\end{gather}
where $C_{1}$ is an integration constant. Let $C_{1}=-F(x_{0})$, where $F(x)$ is a function that has property $F'(x)=f(x)$. Then the equation \eqref{eq1b} can be represented as
$$
(\log |y(x)|)'+\left(\int_{x_0}^x f(t)dt\right)'=\frac{g(x)}{y(x)}, \quad \Rightarrow
$$
$$
(\log |y(x)|)'+\left(\log e^{\int_{x_0}^x f(t)dt}\right)'=\frac{g(x)}{y(x)}, \quad \Rightarrow
$$
$$
\left(\log |y(x)|+\log e^{\int_{x_0}^x f(t)dt}\right)'=\frac{g(x)}{y(x)}, \quad \Rightarrow
$$
$$
\left(\log \left(|y(x)| e^{\int_{x_0}^x f(t)dt}\right)\right)'=\frac{g(x)}{y(x)}, \quad \Rightarrow
$$
$$
\frac{\left(y(x) e^{\int_{x_0}^x f(t)dt}\right)'}{y(x) e^{\int_{x_0}^x f(t)dt}}=\frac{g(x)}{y(x)}, \quad \Rightarrow
$$
$$
\left(y(x) e^{\int_{x_0}^x f(t)dt}\right)' =g(x) e^{\int_{x_0}^x f(\sigma)d\sigma}, \quad \Rightarrow
$$
$$
y(x) e^{\int_{x_0}^x f(t)dt} =\int_{x_0}^x g(\tau) e^{\int_{x_0}^\tau f(\sigma)d\sigma} d\tau + C, \quad \Rightarrow
$$
\begin{gather}\label{eq2c}
y(x)  =e^{-\int_{x_0}^x f(t)dt}\left(\int_{x_0}^x g(\tau) e^{\int_{x_0}^\tau f(\sigma)d\sigma} d\tau + C\right).
\end{gather}
If in the equation \eqref{eq2c} we let $C=y(x_{0})$, then we have the formula \eqref{eq2d}.

{\it 1.2. Bernoulli Differential equation:}
\begin{gather}\label{eq3}
y'+f(x)y=g(x)y^{\alpha},
\end{gather}
where $\alpha\in \mathbb{R}\backslash\{0,1\}$.

{\it Logarithmic integration method.} Let $y$ is not identically zero. Then from the equations \eqref{eq3} we obtain
$$
\frac{y'}{y}+f(x)=\frac{g(x)}{y}y^{\alpha}, \quad \Rightarrow
$$
$$
(\log |y|)'+f(x)=\frac{g(x)}{y}y^{\alpha}, \quad \Rightarrow
$$
$$
(\log |y|)'+\left(\int f(x)dx\right)'=\frac{g(x)}{y}y^{\alpha}, \quad \Rightarrow
$$
$$
(\log |y|)'+\left(\log e^{\int f(x)dx}\right)'=\frac{g(x)}{y}y^{\alpha}, \quad \Rightarrow
$$
$$
\left(\log |y|+\log e^{\int f(x)dx}\right)'=\frac{g(x)}{y}y^{\alpha}, \quad \Rightarrow
$$
\begin{gather}\label{eq3s33434ss}
\left(\log \left(|y| e^{\int f(x)dx}\right)\right)'=\frac{g(x)}{y}y^{\alpha}, \quad \Rightarrow
\end{gather}
$$
\frac{\left(y e^{\int f(x)dx}\right)'}{y e^{\int f(x)dx}}=\frac{g(x)}{y}y^{\alpha}, \quad \Rightarrow
$$
$$
\left(y e^{\int f(x)dx}\right)' =g(x) e^{\int f(x)dx}y^{\alpha}, \quad \Rightarrow
$$
$$
\left(y e^{\int f(x)dx}\right)' =g(x) e^{(1-\alpha)\int f(x)dx}y^{\alpha}e^{\alpha\int f(x)dx}, \quad \Rightarrow
$$
$$
\left(y e^{\int f(x)dx}\right)' =g(x) e^{(1-\alpha)\int f(x)dx}\left(ye^{\int f(x)dx}\right)^{\alpha}, \quad \Rightarrow
$$
$$
\frac{\left(y e^{\int f(x)dx}\right)'}{\left(ye^{\int f(x)dx}\right)^{\alpha}} =g(x) e^{(1-\alpha)\int f(x)dx}, \quad \Rightarrow
$$
$$
\left(\frac{1}{1-\alpha}\left(y e^{\int f(x)dx}\right)^{1-\alpha}\right)' =g(x) e^{(1-\alpha)\int f(x)dx}, \quad \Rightarrow
$$
$$
\frac{1}{1-\alpha}\left(\left(y e^{\int f(x)dx}\right)^{1-\alpha}\right)' =g(x) e^{(1-\alpha)\int f(x)dx}, \quad \Rightarrow
$$
\begin{gather}\label{eq3s33434ss1}
\left(\left(y e^{\int f(x)dx}\right)^{1-\alpha}\right)' =\left(1-\alpha\right)g(x) e^{(1-\alpha)\int f(x)dx}, \quad \Rightarrow
\end{gather}
$$
\left(y e^{\int f(x)dx}\right)^{1-\alpha} =\left(1-\alpha\right)\int g(x) e^{(1-\alpha)\int f(x)dx}dx +C, \quad \Rightarrow
$$
$$
y e^{\int f(x)dx} =\left(\left(1-\alpha\right)\int g(x) e^{(1-\alpha)\int f(x)dx}dx +C\right)^{\frac{1}{1-\alpha}}, \quad \Rightarrow
$$
\begin{gather}\label{eq4}
y =e^{-\int f(x)dx} \left(\left(1-\alpha\right)\int g(x) e^{(1-\alpha)\int f(x)dx}dx +C\right)^{\frac{1}{1-\alpha}}.
\end{gather}

{\it Remark 1.2.1.} At the beginning of the course of the method, we assumed that $y$ be not identically zero $0$. It follows that the equation \eqref{eq3} has a particular solution $y=0$, if $\alpha\in(0,1)$.

{\it Remark 1.2.2.} {\it (The second version of the logarithmic method.) } In the equation \eqref{eq3} we obtain
$$
\frac{y'}{y}+f(x)=\frac{g(x)}{y}y^{\alpha}, \quad \Rightarrow
$$
$$
(\log |y|)'+f(x)=g(x)y^{\alpha-1}, \quad \Rightarrow
$$
$$
(1-\alpha)(\log |y|)'+(1-\alpha)f(x)=(1-\alpha)g(x)y^{\alpha-1}, \quad \Rightarrow
$$
$$
((1-\alpha)\log |y|)'+(1-\alpha)f(x)=(1-\alpha)g(x)y^{\alpha-1}, \quad \Rightarrow
$$
$$
(\log |y|^{1-\alpha})'+(1-\alpha)f(x)=(1-\alpha)g(x)y^{\alpha-1}, \quad \Rightarrow
$$
$$
\frac{ \left(y^{1-\alpha}\right)' }{ y^{1-\alpha} }   +(1-\alpha)f(x)=(1-\alpha)g(x)y^{\alpha-1}, \quad \Rightarrow
$$
$$
\left(y^{1-\alpha}\right)'  +(1-\alpha)f(x)y^{1-\alpha}=(1-\alpha)g(x)y^{\alpha-1}y^{1-\alpha}, \quad \Rightarrow
$$
\begin{gather}\label{eq4s1112}
\left(y^{1-\alpha}\right)'  +(1-\alpha)f(x)y^{1-\alpha}=(1-\alpha)g(x).
\end{gather}
The equation \eqref{eq4s1112} is a linear inhomogeneous first order differential equation, with respect to the function $y^{1-\alpha}$. Its solution by the with formula \eqref{eq2}, has the form
\begin{gather}\label{eq4s1114}
y^{1-\alpha} =e^{-(1-\alpha)\int f(x)dx} \left(\left(1-\alpha\right)\int g(x) e^{(1-\alpha)\int f(x)dx}dx +C\right).
\end{gather}
The formula \eqref{eq4s1114} implies the solution \eqref{eq4}.

{\it Remark 1.2.3.} {\it (The third version of the logarithmic method.) } In the equation \eqref{eq3s33434ss} we obtain

$$
(1-\alpha)\left(\log \left(|y| e^{\int f(x)dx}\right)\right)'=(1-\alpha)g(x)y^{\alpha-1}, \quad \Rightarrow
$$
$$
\left((1-\alpha)\log \left(|y| e^{\int f(x)dx}\right)\right)'=(1-\alpha)g(x)y^{\alpha-1}, \quad \Rightarrow
$$
$$
\left(\log \left(|y| e^{\int f(x)dx}\right)^{1-\alpha}\right)'=(1-\alpha)g(x)y^{\alpha-1}, \quad \Rightarrow
$$
$$
\frac{\left(\left(y e^{\int f(x)dx}\right)^{1-\alpha }\right)'}{\left(y e^{\int f(x)dx}\right)^{1-\alpha }}=(1-\alpha)g(x)y^{\alpha-1}, \quad \Rightarrow
$$
\begin{gather}\label{eq4s1114dd1}
\left(\left(y e^{\int f(x)dx}\right)^{1-\alpha }\right)'=(1-\alpha)g(x)y^{\alpha-1}\left(y e^{\int f(x)dx}\right)^{1-\alpha }=(1-\alpha)g(x)e^{(1-\alpha)\int f(x)dx}.
\end{gather}
The equation \eqref{eq4s1114dd1} is similar to the equation \eqref{eq3s33434ss1}.

{\it 1.3. The equation of the form:}
\begin{gather}\label{eq5}
y'+f(x)e^{\beta y}=g(x),
\end{gather}
where $\beta\in \mathbb{R}\backslash{\{0\}}$.

{\it Logarithmic integration method.} In the equation \eqref{eq5} we get
$$
\left(\log\left(e^{y}\right)\right)'+f(x)e^{\beta y}=g(x), \Rightarrow
$$
$$
-\beta\left(\log\left(e^{y}\right)\right)' -\beta f(x)e^{\beta y}=-\beta g(x), \Rightarrow
$$
$$
\left(-\beta\log\left(e^{y}\right)\right)' -\beta f(x)e^{\beta y}=-\beta g(x), \Rightarrow
$$
$$
\left(\log\left(e^{-\beta y}\right)\right)' -\beta f(x)e^{\beta y}=-\beta g(x), \Rightarrow
$$
$$
\frac{\left(e^{-\beta y}\right)'}{e^{-\beta y}} -\beta f(x)e^{\beta y}=-\beta g(x), \Rightarrow
$$
$$
\left(e^{-\beta y}\right)' -\beta f(x)e^{\beta y}e^{-\beta y}=-\beta g(x)e^{-\beta y}, \Rightarrow
$$
$$
\left(e^{-\beta y}\right)' -\beta f(x)=-\beta g(x)e^{-\beta y}, \Rightarrow
$$
\begin{gather}\label{eq6}
\left(e^{-\beta y}\right)' +\beta g(x)e^{-\beta y}=\beta f(x).
\end{gather}
The equation \eqref{eq6} is a linear inhomogeneous first order differential equation, with respect to the function $e^{-\beta y}$. Its solution, by the formula \eqref{eq2}, has the form
\begin{gather}\label{eq7}
e^{-\beta y}  =e^{-\beta\int g(x)dx} \left( \beta\int f(x) e^{\beta\int g(x)dx} dx + C \right).
\end{gather}
Solving the equation \eqref{eq7}, with respect to $y$, we have
$$
y  =-\frac{1}{\beta}\log\left(e^{-\beta\int g(x)dx} \left( \beta\int f(x) e^{\beta\int g(x)dx} dx + C \right)\right), \Rightarrow
$$
\begin{gather}\label{eq8}
y  =\int g(x)dx-\frac{1}{\beta}\log \left( \beta\int f(x) e^{\beta\int g(x)dx} dx + C \right).
\end{gather}

\medskip
\section{Second order differential equations}

{\it 2.1. Linear homogeneous second order differential equation:}
\begin{gather}\label{eq9}
y''+by'+cy=0,
\end{gather}
where $b\in \mathbb{R}$, $c\in \mathbb{R}$.

Let $y$ is not identically zero. Then from the equation \eqref{eq9} we obtain

$$
\frac{y''}{y}+b\frac{y'}{y}+c=0, \quad \Rightarrow
$$
\begin{gather}\label{eq10}
(\log |y|)''+((\log |y|)')^2+ b(\log |y|)'+c=0,
\end{gather}
because $\frac{y'}{y}=(\log |y|)'$, $(\log |y|)''=\left(\frac{y'}{y}\right)'=\frac{y''}{y}-\left(\frac{y'}{y}\right)^2=\frac{y''}{y}-\left((\log |y|)'\right)^2$, $\Rightarrow$ $\frac{y''}{y}=(\log |y|)''+\left((\log |y|)'\right)^2$.
Let in the equation \eqref{eq10}:
$$
(\log |y|)'=z.
$$
Then we have equation \eqref{eq10} in the form
\begin{gather}\label{eq11}
z'+z^2+bz+c=0, \quad \Rightarrow
\end{gather}
\begin{gather}\label{eq12}
\frac{z'}{z^2+bz+c}=-1.
\end{gather}

{\it Case 1.} $b^2-4c>0$. In this case we have equation \eqref{eq12} has be form
$$
\frac{z'}{\left(z+\frac{b}{2}\right)^2-\frac{1}{4}(b^2-4c)}=-1, \quad \Rightarrow
$$
$$
\frac{z'}{\left(z+\frac{b}{2}\right)^2-\left(\frac{1}{2}\sqrt{b^2-4c}\right)^2}=-1, \quad \Rightarrow
$$
\begin{gather}\label{eq13}
\frac{1}{2\frac{1}{2}\sqrt{b^2-4c}}\log\left|\frac{z+\frac{b}{2}-\frac{1}{2}\sqrt{b^2-4c}}{z+\frac{b}{2}+\frac{1}{2}\sqrt{b^2-4c}}\right|=-x+C_{1,1}.
\end{gather}
Let in the equation \eqref{eq13}: $z+\frac{b}{2}=\xi$, $\frac{1}{2}\sqrt{b^2-4c}=\gamma$. Then we obtain
$$
\frac{1}{\sqrt{b^2-4c}}\log\left|\frac{\xi-\gamma}{\xi+\gamma}\right|=-x+C_{1,1}, \quad \Rightarrow
$$
$$
\frac{\xi-\gamma}{\xi+\gamma}=C_{1,2}e^{-\sqrt{b^2-4c}x}, \quad C_{1,2}=e^{C_{1,1}}, \quad \Rightarrow
$$
$$
1-\frac{2\gamma}{\xi+\gamma}=C_{1,2}e^{-\sqrt{b^2-4c}x}, \quad \Rightarrow
$$
$$
\frac{1}{\xi+\gamma}=\frac{1}{2\gamma}+C_{1,3}e^{-\sqrt{b^2-4c}x}, \quad C_{1,3}=-\frac{1}{2\gamma}C_{1,2}, \quad \Rightarrow
$$
\begin{gather}\label{eq14}
\xi+\gamma=\frac{1}{\frac{1}{2\gamma}+C_{1,3}e^{-\sqrt{b^2-4c}x}}.
\end{gather}
Returning to the change of variables $z+\frac{b}{2}=\xi$, $\frac{1}{2}\sqrt{b^2-4c}=\gamma$ in the equation \eqref{eq14}, we obtain
$$
z+\frac{b}{2}+\frac{1}{2}\sqrt{b^2-4c}=\frac{1}{\frac{1}{\sqrt{b^2-4c}}+C_{1,3}e^{-\sqrt{b^2-4c}x}}, \quad \Rightarrow
$$
$$
z=-\left(\frac{b}{2}+\frac{1}{2}\sqrt{b^2-4c}\right)+\frac{1}{\frac{1}{\sqrt{b^2-4c}}+C_{1,3}e^{-\sqrt{b^2-4c}x}}, \quad \Rightarrow
$$
$$
z=-\left(\frac{b}{2}+\frac{1}{2}\sqrt{b^2-4c}\right)+\frac{\sqrt{b^2-4c}}{1+C_{1,4}e^{-\sqrt{b^2-4c}x}}, \quad C_{1,4}=\sqrt{b^2-4c}C_{1,3}, \quad \Rightarrow
$$
\begin{gather}\label{eq15}
z=-\left(\frac{b}{2}+\frac{1}{2}\sqrt{b^2-4c}\right)+\frac{e^{\sqrt{b^2-4c}x}\sqrt{b^2-4c}}{e^{\sqrt{b^2-4c}x}+C_{1,4}}.
\end{gather}
Because $z=(\log |y|)'$, then we have in the equation \eqref{eq15}
$$
(\log |y|)'=-\left(\frac{b}{2}+\frac{1}{2}\sqrt{b^2-4c}\right)+\frac{e^{\sqrt{b^2-4c}x}\sqrt{b^2-4c}}{e^{\sqrt{b^2-4c}x}+C_{1,4}}, \quad \Rightarrow
$$
$$
\log |y|=\int\frac{e^{\sqrt{b^2-4c}x}\sqrt{b^2-4c}}{e^{\sqrt{b^2-4c}x}+C_{1,4}}dx-\left(\frac{b}{2}+\frac{1}{2}\sqrt{b^2-4c}\right)x+\log|C_{2,1}|=
$$
$$
=\log \left|e^{\sqrt{b^2-4c}x}+C_{1,4}\right|-\left(\frac{b}{2}+\frac{1}{2}\sqrt{b^2-4c}\right)x+\log|C_{2,1}|, \quad \Rightarrow
$$
$$
y=\left(e^{\sqrt{b^2-4c}x}+C_{1,4}\right)e^{-\left(\frac{b}{2}+\frac{1}{2}\sqrt{b^2-4c}\right)x}C_{2,1}, \quad \Rightarrow
$$
\begin{gather}\label{eqd141}
y=C_{1}e^{\left(-\frac{b}{2}-\frac{1}{2}\sqrt{b^2-4c}\right)x}+C_{2}e^{\left(-\frac{b}{2}+\frac{1}{2}\sqrt{b^2-4c}\right)x},
\end{gather}
where $C_{1}=C_{2,1}C_{1,4}$, $C_{2}=C_{2,1}$ is an integration constant.

{\it Case 2.} $b^2-4c=0$.
In this case we have equation \eqref{eq12} has be form
$$
\frac{z'}{\left(z+\frac{b}{2}\right)^2}=-1.
$$
Step by step from the last equation we obtain
$$
\left(-\frac{1}{z+\frac{b}{2}}\right)'=-1, \quad \Rightarrow
$$
$$
\frac{1}{z+\frac{b}{2}}=x+C_{1,5}, \quad \Rightarrow
$$
$$
z+\frac{b}{2}=\frac{1}{x+C_{1,5}}, \quad \Rightarrow
$$
\begin{gather}\label{eq15sd1}
z=-\frac{b}{2}+\frac{1}{x+C_{1,5}}.
\end{gather}
Because $z=(\log |y|)'$, then we have in the equation \eqref{eq15sd1}
$$
(\log |y|)'=-\frac{b}{2}+\frac{1}{x+C_{1,5}}, \quad \Rightarrow
$$
$$
\log |y|=-\frac{b}{2}x+\log \left|x+C_{1,5}\right|+\log\left|C_{2}\right|, \quad \Rightarrow
$$
\begin{gather}\label{eqd142}
y=e^{-\frac{b}{2}x}\left(x+C_{1,5}\right)C_{2}=C_{1}e^{-\frac{b}{2}x}+C_{2}xe^{-\frac{b}{2}x},
\end{gather}
where $C_{1}=C_{1,5}C_{2}$ is an integration constant.

{\it Case 3.} $b^2-4c<0$. In this case we have equation \eqref{eq12} has be form
$$
\frac{z'}{\left(z+\frac{b}{2}\right)^2+\frac{1}{4}(4c-b^2)}=-1, \quad \Rightarrow
$$
$$
\frac{z'}{\left(z+\frac{b}{2}\right)^2+\left(\frac{1}{2}\sqrt{4c-b^2}\right)^2}=-1, \quad \Rightarrow
$$
$$
\frac{1}{\frac{1}{2}\sqrt{4c-b^2}}\arctan \frac{z+\frac{b}{2}}{\frac{1}{2}\sqrt{4c-b^2}}=-x+C_{1,6} \quad \Rightarrow
$$
$$
\arctan \frac{z+\frac{b}{2}}{\frac{1}{2}\sqrt{4c-b^2}}=-\frac{1}{2}\sqrt{4c-b^2}x+C_{1,7}, \quad C_{1,7}=C_{1,6}\frac{1}{2}\sqrt{4c-b^2}, \quad \Rightarrow
$$
$$
\frac{z+\frac{b}{2}}{\frac{1}{2}\sqrt{4c-b^2}}=\tan \left(-\frac{1}{2}\sqrt{4c-b^2}x+C_{1,7}\right), \quad \Rightarrow
$$
\begin{gather}\label{eq15sd2}
z=-\frac{b}{2}+\frac{1}{2}\sqrt{4c-b^2}\tan \left(-\frac{1}{2}\sqrt{4c-b^2}x+C_{1,7}\right).
\end{gather}
Because $z=(\log |y|)'$, then we have in the equation \eqref{eq15sd2}
$$
(\log |y|)'=-\frac{b}{2}+\frac{1}{2}\sqrt{4c-b^2}\tan \left(-\frac{1}{2}\sqrt{4c-b^2}x+C_{1,7}\right), \quad \Rightarrow
$$
$$
\log |y|=-\frac{b}{2}x+\frac{1}{2}\sqrt{4c-b^2}\int\tan \left(-\frac{1}{2}\sqrt{4c-b^2}x+C_{1,7}\right)dx+\log\left|C_{2,2}\right|=
$$
$$
=-\frac{b}{2}x+\log\left|\cos\left(-\frac{1}{2}\sqrt{4c-b^2}x+C_{1,7}\right)\right|+\log\left|C_{2,2}\right|, \quad \Rightarrow
$$
$$
y=e^{-\frac{b}{2}x}\cos\left(-\frac{1}{2}\sqrt{4c-b^2}x+C_{1,7}\right)C_{2,2}=
$$
$$
=C_{2,2}e^{-\frac{b}{2}x}\left(\cos\left(-\frac{1}{2}\sqrt{4c-b^2}x\right)\cos \left(C_{1,7}\right)-\sin\left(-\frac{1}{2}\sqrt{4c-b^2}x\right)\sin \left(C_{1,7}\right)\right)=
$$
\begin{gather}\label{eqd143}
=e^{-\frac{b}{2}x}\left(C_{1}\cos\left(\frac{1}{2}\sqrt{4c-b^2}x\right)+C_{2}\sin\left(\frac{1}{2}\sqrt{4c-b^2}x\right)\right),
\end{gather}
where $C_{1}=C_{2,2}\cos (C_{1,7})$, $C_{2}=C_{2,2}\sin (C_{1,7})$ is an integration constant.

The formulas \eqref{eqd141}, \eqref{eqd142}, \eqref{eqd143}, solve the equation \eqref{eq9} in the respective cases 1,2,3.
\medskip

\textbf{Conclusion.} This method in chapter 2 makes it possible to obtain these solutions without applying a complex analysis and finding a solution in the form $y=\psi(x)e^{\zeta x}$. Also, we got exact solutions for many kinds of first-order differential equations in chapter 1.

\end{document}